\documentclass[10pt,reqno,twoside]{amsart}

 \renewcommand{\d}{{\mathrm{d}}} %differential
 \renewcommand{\i}{{\mathrm{i}}} %\sqrt{-1}

\renewcommand{\Re}{{\mathfrak{Re}}}
\renewcommand{\Im}{{\mathfrak{Im}}}

\newcommand{\E}{\mathbb{E}}
\newcommand{\Prob}{\mathbb{P}}
\DeclareMathOperator{\meas}{meas}
\DeclareMathOperator{\rank}{rank}

\def\({\left(}
\def\){\right)}

\newtheorem{theorem}{Theorem}[section]

\newtheorem{lemma}[theorem]{Lemma}
\newtheorem{conj}{Conjecture}

 \numberwithin{equation}{section}
 \numberwithin{theorem}{section}

\begin{document}

\title{The maximum size of $L$-functions}

\author{David W. Farmer}
 \address{American Institute of Mathematics, 360 Portage Avenue, Palo Alto,
CA 94306-2244, USA.}
 \email{farmer@aimath.org}

\author{S.M. Gonek}
\address{Department of Mathematics, University of Rochester, Rochester, NY 14627, USA}
 \email{gonek@math.rochester.edu}

\author{C.P. Hughes}
 \address{Department of Mathematics, University of
Michigan, Ann Arbor, MI 48109-1043, USA}
 \curraddr{Department of Mathematics, University of York, Heslington, York, YO10 5DD, U.K.}
 \email{hughes@aimath.org}

\thanks{All three authors were supported by an NSF Focused Research
Group grant DMS 0244660. Work of the second author was also
supported by NSF grant DMS 0201457. Work of the third author was
partially supported by EPSRC grant N09176. All three authors wish
to thank the Isaac Newton Institute for their hospitality during
the course of this work.}

  \dedicatory{Dedicated to Hugh Montgomery on his 60th birthday.}

  \date{August 1, 2006}

\begin{abstract}

We conjecture the true rate of growth of the maximum size of the
Riemann zeta-function and other $L$-functions. We support our
conjecture using arguments from random matrix
theory, conjectures for moments of $L$-functions,
and also by assuming a random model for the primes.

\end{abstract}

\maketitle

\section{Introduction and statement of results}\label{sect1}

A fundamental problem in analytic number theory is to calculate the
maximum size of $L$-functions in the critical strip. For example,
the importance of the Lindel\"of Hypothesis, which is a consequence of the
Riemann Hypothesis, is that it provides at least a crude estimate
for the maximum in the case of the Riemann zeta-function.
In this paper we  use a variety of methods to
conjecture the true rate of growth.

Consider first the Riemann zeta-function, which is a prototypical
$L$-function. The Lindel\"of Hypothesis  asserts that for every
$\varepsilon>0$, $\zeta(\tfrac12+\i t) = O(t^\varepsilon)$ (here
we assume  $t$ is positive). Under the Riemann Hypothesis, one can
show that
\begin{equation}\label{eq:zeta_bigO}
\zeta(\tfrac12+\i t) = O\left(\exp\left(C\frac{\log t}{\log \log t}\right)\right)
\end{equation}
for some constant $C$ (see Theorem 14.14A of \cite{Ti}, for
example).
Several results of the form
\begin{equation}\label{eq:zeta_Omega}
|\zeta(\tfrac12+\i t)|= \Omega\left(\exp\left(C'\sqrt{\mathstrut \frac{\log t}{\log \log
t}}\right)\right)
\end{equation}
have also been established.
Assuming the Riemann Hypothesis (RH), Montgomery~\cite{Mon} showed $C'\ge 1/20$.
Balasubramanian and
Ramachandra~\cite{BalRam} improved the constant ~$C'$ and removed the
assumption of RH.
Soundararajan~\cite{S} further improved the estimate to~$C'\ge 1$,
and he has also obtained similar results for the central values
of a family of $L$-functions.
In fact,
Soundararajan's calculations show that the proportion of $t$ for
which $\zeta(\tfrac12+\i t)$ is this big is quite large,
suggesting that it may get bigger still.
Numerical calculations of Kotnik~\cite{K} indicate that $C'>2$
and perhaps $C'$ can be much larger.

We are interested in finding out which of equations
\eqref{eq:zeta_bigO} or \eqref{eq:zeta_Omega} is closer to the
truth. This is part of a class of problems that has recently come to be
known as the ``1 or 2?'' question, where one has an $O$-result and an $\Omega$-result
which, suitably interpreted, differ by a factor of~$2$.  In the
case here, the unknown factor is the power of $\log t$ in the
exponential. The calculations in this paper support the view that
``1'' is the correct answer in this case, and we make the
following conjecture:
\begin{conj}\label{conj:max value zeta}
\begin{equation}
\max_{t \in [0,T]} |\zeta(\tfrac12+\i t)| =
\exp\left((1+o(1))\sqrt{\mathstrut \tfrac12 \log T \log\log T}
\right) .
\end{equation}
\end{conj}

Similar arguments to those presented for $\zeta(\frac12+\i t)$
work for $S(t)$, the error term in the number of zeros of the
zeta-function up to height $t$, and lead to
\begin{conj} \label{conj:max value S}
\begin{equation}
\limsup_{t\to\infty} \frac{S(t)}{\sqrt{\mathstrut \mathstrut\log t \log\log t}} =
\frac{1}{\pi\sqrt{\mathstrut 2}} .
\end{equation}
\end{conj}

Much recent progress in understanding analytic properties of
$L$-functions has come from the idea of a ``family'' of
$L$-functions with an associated symmetry type.
The idea is that to a collection of $L$-functions, with
appropriate natural conditions, one can associate a
classical compact group: unitary, symplectic, or orthogonal.
One expects the analytic properties of the $L$-functions
to be largely governed only by the symmetry type.
Here we apply this philosophy to conjecture the maximal size
of the critical values of $L$-functions.

A family $\mathcal F$ of $L$-functions is partially ordered by the
``conductor'' $c(F)$ for $F\in \mathcal F$.  Our calculations
assume that $\#\{F\in \mathcal F \ :\ c(F)<D\} \approx D$.
Straightforward modifications can handle the case in which the
family grows like $D^A$ for any $A>0$.

For a more detailed discussion of families of $L$-functions,
see~\cite{CF, CFKRS}.  (However, note that~\cite{CFKRS} introduces
a refined notion of ``conductor'',  which, asymptotically, is
the logarithm of the ``usual'' conductor  used here.)

\begin{conj}\label{conj:max value L}
Suppose $\mathcal F$ is a family of $L$-functions  and, for $F\in
\mathcal F$, let $c(F)$ denote the conductor of~$F$. With $B=1/2$
for unitary families and $B=1$ for symplectic and orthogonal
families, we have
\begin{equation}
\max_{\substack{F\in \mathcal F \\ c(F) \le D}} |F(\tfrac12)| =
\exp\left((1+o(1)) \sqrt{\mathstrut B \log D \log\log D }  \right)
.
\end{equation}
The implied constant depends only on $\mathcal F$.
\end{conj}

For example, for the symplectic family of real primitive
Dirichlet $L$-functions, $L(s,\chi_d)$, where $\chi_d=\binom{d}{\cdot}$,
we conjecture that
\begin{equation}
\max_{|d| \le D} |L(\tfrac12 ,\chi_d)| =
\exp\left(
      (1+o(1))  \sqrt{\mathstrut \log D \log\log D } \right) .
\end{equation}

Similarly, for the orthogonal family of Dirichlet series
associated with holomorphic cusp forms, $L(s,f)$, where $f\in
S_k(\Gamma_0(N))$, the conductor is $kN$, so we conjecture that
\begin{equation}
\max_{\substack{f\in S_k(\Gamma_0(N)) \\ kN \le D}} |L(\tfrac{k}2 ,f)| =
\exp\left(
      (1+o(1)) \sqrt{\mathstrut \log  D  \log\log  D  }  \right).
\end{equation}

Note that Conjecture~\ref{conj:max value L} contains
Conjecture~\ref{conj:max value zeta}, because any primitive
$L$-function, $L(s)$, has associated with it the unitary family
\begin{equation}
\mathcal F_L := \left\{ F_y(s):=L(s+i y) \ \ |\  \ y\in\mathbb R \right\}
\end{equation}
with conductor~$c(F_y)\sim |y|$.

Our conjectures suggest that on the critical line the answer to
the ``1 or 2'' question is ``1''.  Work of
Montgomery and Vaughan~\cite{MV}
and Granville and Soundararajan~\cite{GS1, GS2} has suggested that
the answer also is ``1'' on the $1$-line.
Thus in both cases, the maximum
value the $L$-function attains appears to be closer to the
$\Omega$-result than to the
$O$-result.

In Section~\ref{sect:new model approach} we use a rigorous
approximation to the zeta-function due to Gonek, Hughes, and
Keating~\cite{GHK} to justify Conjecture~{\ref{conj:max value zeta}.
This approximation represents  $\zeta(s)$ as a product over
primes times a product over zeros. We use characteristic
polynomials of random unitary matrices to model the product over
zeros, and a separate probabilistic model due to Granville and
Soundararajan for the product over primes. The approximation to
the zeta-function has a parameter which controls the relative
contribution of the primes and the zeros.  We show that the
predicted maximal order of the zeta-function is the same
independent of the choice of parameter.  That is, whether we use
only the primes, or only the zeros, or some combination of the two, we obtain
Conjecture~{\ref{conj:max value zeta}.

In Section~\ref{sect:moments} we use recent conjectures for
moments of $L$-functions to give an alternative justification for
Conjecture~\ref{conj:max value zeta}.  Our approach also provides
new limits on the range of validity of those conjectured moments.

In Section~\ref{sect:S} we modify the treatment in
Section~\ref{sect:new model approach} to obtain
Conjecture~\ref{conj:max value S}.

In Section~\ref{sect:other} we describe how to extend our approach
to obtain Conjecture~\ref{conj:max value L}. We also describe some
other approaches to obtaining these conjectures and  then indicate
possible arguments against the conjectures.

Finally, in Appendix~\ref{sect:proofs} we prove a theorem about
random matrix polynomials that is used in Section~\ref{sect:new
model approach}.

We thank Andrew Granville and Soundararajan for allowing us to
incorporate their work on extreme values using the primes, and we
thank Andrew Booker, Brian Conrey, Hugh Montgomery and Doug Ulmer for helpful
conversations.

%%%%%%%%%%%  !!!!!!!!!!!  %%%%%%%%%%%%%%%
%%%%%%%%%%%  New Section  %%%%%%%%%%%%%%%
%%%%%%%%%%%  !!!!!!!!!!!  %%%%%%%%%%%%%%%

\section{A probabilistic model for the zeta-function}\label{sect:new model approach}

Gonek, Hughes, and Keating~\cite{GHK} have proved that if $s=\sigma+\i
t$, with $0\leq\sigma\leq 1$ and $|t|\geq 2$, then for $X>2$ and $K$
any positive integer,
\begin{equation}\label{eq:zeta= PZ}
    \zeta(s) = P_X(s) Z_X(s) \left(1+O\left(\frac{X^{2-\sigma+K}}{(|t|
    \log X)^{K}}\right)+O(X^{-\sigma}\log X)\right)\;,
\end{equation}
where
\begin{equation}\label{eq:P_X}
P_X(s) := \exp\left(\sum_{n\leq X}\frac{\Lambda(n)}{n^{s}\log n}
\right)\;,
\end{equation}
$\Lambda(n)$ is the von-Mangoldt function, and
\begin{equation}\label{eq:Z_X}
Z_X(s):= \exp\left(-\sum_{\rho}U((s-\rho)\log X)\right)\;.
\end{equation}
Here the $\rho$ are non-trivial zeros of $\zeta(s)$ and $U(z) =
\int_0^\infty u(x) E_1(z\log x)\;\d x$, where $E_1(z) =
\int_{z}^{\infty} \frac{e^{-w}}{w}\;\d w$ is the exponential
integral and $u$ is any smooth function supported in $[e^{1-1/X}
, e]$.

The parameter $X$ controls the relative influence of the primes
and the zeros. If $X$ is large, there are many primes in $P_X(s)$,
and only the zeros very close to $s$ effect the product in
$Z_X(s)$, while if $X$ is small, the zeros further away from $s$
make a contribution to $Z_X(s)$, but the number of primes in
$P_X(s)$ is diminished. When $X$ is not too large,
we expect  $Z_X$ and $P_X$ to behave somewhat
independently,  and  Gonek, Hughes, and Keating~\cite{GHK}
give evidence of this. In the remainder of this section we describe
probabilistic models for $P_X$ and $Z_X$ which, assuming independence, will give
Conjecture~\ref{conj:max value zeta}.
In Section~\ref{sect:char poly} we describe our model
for the large values of~$|Z_X|$,
establish some new results on the size of characteristic
polynomials of random unitary matrices, and justify
Conjecture~\ref{conj:max value zeta} by choosing $X$ small.  In
Section~\ref{sect:prime model} we describe Granville and Soundararajan's
model for $P_X$ and justify Conjecture~\ref{conj:max value zeta} by choosing
$X$ large.  Then, in Section~\ref{sect:primes and zeros}
we combine $Z_X$ and $P_X$, showing that intermediate values of $X$
also lead to Conjecture~\ref{conj:max value zeta}.

\subsection{A random matrix model for large values of $Z_X(\frac12+\i t)$}\label{sect:char poly}

Here we study the characteristic polynomial
\begin{equation}
\Lambda_U(\theta) = \det\left(I-U e^{-\i\theta}\right) =
\prod_{n=1}^N \left(1-e^{\i(\theta_n-\theta)}\right)
\end{equation}
of a random unitary matrix $U\in U(N)$  chosen uniformly with respect to Haar
measure. The characteristic polynomial $\Lambda_U(\theta)$ was
first developed as a model for the Riemann zeta-function by
Keating and Snaith~\cite{KeaSna}. In \cite{GHK}, building
on~\cite{KeaSna}, it is argued that for $t \approx T$, $Z_X(\frac12+\i t)$, given by
\eqref{eq:Z_X}, can be modeled by $\Lambda_U(\theta)$, where $U\in
U(N)$ with
\begin{equation}\label{eqn:N}
N = \left[ \frac{\log T }{ e^\gamma \log X} \right].
\end{equation}
We will prove a result about the value distribution and maximal
size of~$|\Lambda_U(\theta)|$ in Appendix ~\ref{sect:proofs}, and use
it to conjecture the distribution of large values of~$|Z_X|$.

The largest value of $|\Lambda_U(\theta)|$  is $2^N$, and values near
this  occur  when $U$ is in a small neighborhood of scalar multiples of
the identity matrix.
If $X = e^{o(\log\log T)}$, this violates the known bound  on
$|\zeta(\frac12+\i t)|$, so our model for the large values of
$|Z_X(\frac12+\i t)|$ must do something more subtle than just take
the maximum of $|\Lambda_U|$ over all $U\in U(N)$.

If $T$ and $X$ are
thought of as fixed, then matrices of size
$N = [\frac{\log T}{e^\gamma \log X}]$
should model the zeta-function as long as
$T /X^{e^\gamma}  < t <T$.
If $X>2$, say, then up to constants there are $T\log T$
zeros in this interval. Therefore, in order to
have the same number of eigenvalues, one needs
\begin{equation}\label{eqn:M}
M=[T \log T / N] \approx \exp\left({e^\gamma N \log X}\right) \log X
\end{equation}
matrices. Thus, one
plausible guess for the maximum value of $|Z_X(\frac12+\i t)|$ for
$0<t<T$  is $K = K(M,N)$, where $N$ and $M$ are given
in~\eqref{eqn:N} and ~\eqref{eqn:M}, respectively, and $K$ is
the smallest  possible function of $M$ and $N$ such that
\begin{equation}\label{tail prob}
\Prob\left\{ \max_{1\leq j \leq M}
\max_{\theta}|\Lambda_{U_j}(\theta)| \leq K \right\} \to 1
\end{equation}
as $N\to\infty$. Such a $K$ is found in
Theorem~\ref{thm:max size max Lambda}.

We have glossed over some issues here, but we argue that they are
not essential. First, we have claimed that matrices of size $N$
model $Z_X(\frac12+ \i t)$ for   $T /X^{e^\gamma} < t< T$,
whereas we want $0<t<T$. However, if $X\to\infty$,
then $[T /X^{e^\gamma}, T]$ will
cover almost all of $[0,T]$, and so should capture
the maximum. Secondly, we have been slightly cavalier in dropping
the condition that $N$ should be an integer, which will have an
effect on the number of matrices, $M$, we maximize over. However,  as we
will see below, the answer depends only on the logarithmic size of
$M$, so this is not a serious problem.  Finally, we remark that  the
placement of $e^\gamma$ in our definitions of $N$ and $M$  is actually
irrelevant: our heuristics are sufficiently robust that
increasing $N$ by any fixed constant and decreasing $M$
correspondingly leads to the same conjectured maximum.  We include
the $e^\gamma$ factor to be consistent with~\cite{GHK}, where
the precise choice of $N$ does matter.

We now find an explicit $K$ satisfying \eqref{tail prob}.
\begin{theorem}\label{thm:max size max Lambda}
Fix $\delta>0$. Let $M=\exp(N^\beta)$, with
$\delta<\beta<2-\delta$, and set
\begin{equation}\label{eq:K}
K_\varepsilon(N)=\exp\left(
\left(\sqrt{\mathstrut 1-\tfrac12\beta}+\varepsilon\right) \sqrt{\mathstrut \log M}
\sqrt{\mathstrut \log N} \right) .
\end{equation}
If $U_1,\dots,U_M$ are chosen independently from $U(N)$, then as
$N\to\infty$,
\begin{equation}
\Prob\left\{ \max_{1\leq j \leq M} \max_{\theta}
|\Lambda_{U_j}(\theta)| \leq K_\varepsilon(N) \right\} \to 1
\end{equation}
for all $\varepsilon>0$ and for no $\varepsilon<0$.
\end{theorem}

\begin{proof}
Note that by the independence of the $U_j$,
\begin{equation}
\Prob\left\{ \max_{1\leq j \leq M} \max_{\theta}
|\Lambda_{U_j}(\theta)| \leq K_\varepsilon(N) \right\}
=\Prob\left\{ \max_{\theta} |\Lambda_{U}(\theta)| \leq
K_\varepsilon(N) \right\}^M
\end{equation}
and, for this to tend to $1$ as $N\to\infty$, we must have
\begin{equation}
M \log \left(1-\Prob\left\{ \max_{\theta} |\Lambda_{U}(\theta)|
> K_\varepsilon(N) \right\} \right) \to 0 \ .
\end{equation}
Thus, the proof of the theorem (and all  similar ones in
this paper) requires knowledge of the tails of the distribution,
and this is given by Lemma~\ref{lemma:maxthetadist}.
If $M=\exp(N^\beta)$ and $K=K_\varepsilon(N)$ is as in \eqref{eq:K}
with $0<\beta<2$ and $\varepsilon>-\sqrt{\mathstrut 1-\beta/2}$, then by
Lemma~\ref{lemma:maxthetadist} one  easily finds that
\begin{align}
M \log \left(1-\Prob\left\{ \max_{\theta} |\Lambda_{U}(\theta)|
> K_\varepsilon(N) \right\} \right) &= \exp(N^\beta) \exp\left(-\frac{\log^2 K}{ \log N - \log \log K}
(1+o(1))\right) \\
&\to 0 \notag
\end{align}
as $N\to\infty$ for all $\varepsilon>0$, but for no
$\varepsilon<0$.
\end{proof}

To summarize, we use the characteristic polynomials
$\Lambda_{U}(\theta)$ of random unitary matrices $U\in U(N)$ to
model $Z_X(\frac12+\i t)$. To model the large values of
$|Z_X(\frac12+\i t)|$ for $t\in [0,T]$ we choose $N$ as
in~\eqref{eqn:N} with $\log X < (\log N)^A$ for some $A$, and
we take about
\begin{equation}\label{eqn:Nc}
M= N^c \exp\left(e^\gamma {N \log X}\right)
\end{equation}
different matrices. Here $c\ge 0$ is fixed, and we include it to
allay concerns that choosing too few matrices  may miss some large
values. With these values for $M$ and  $N$,  it follows that  $\beta$
 in  Theorem~\ref{thm:max size max Lambda} is  $\sim 1$, and this leads to the following
conjecture:
\begin{conj}\label{conj:max value ZX}
If $2<X<\log^A T$, then
\begin{equation}
\max_{t\in[0,T]} |Z_X(\tfrac12+\i t)|  = \exp\left( (1+o(1))
\sqrt{\mathstrut \tfrac12 \log T\log \log T}\right)
\end{equation}
as $T\to\infty$.
\end{conj}

We can now complete our argument for  Conjecture~\ref{conj:max value zeta}.
\begin{proof}[Justification of Conjecture~\ref{conj:max value zeta}: ]
By the prime number theorem and \eqref{eq:P_X}, we see that
\begin{equation}
|P_X(\tfrac12+\i t)| \leq \exp\left(\sum_{n\leq X}
\frac{\Lambda(n)}{\sqrt{\mathstrut n} \log n} \right)=
O\left(\exp\left(3 \frac{ \sqrt{\mathstrut X}}{\log X}\right)\right) .
\end{equation}
Thus, if  $X=O(\log T)$ and $T /X^{e^\gamma} < t< T$, then
\begin{equation}
P_X(\tfrac12+\i t)=O\left(
        \exp\left(C\frac{\sqrt{\mathstrut \log t}}{\log \log t}\right)
        \right).
\end{equation}
Combining this,  \eqref{eq:zeta= PZ}, and Conjecture~\ref{conj:max value ZX},
we obtain Conjecture~\ref{conj:max value zeta}.
\end{proof}

The argument above essentially splits the critical line into blocks
of size $1$, maximizes over each block, and then finds the maximum
of the maxima. However, one might instead wish to sample the critical
line at many evenly spaced points. If they are not too sparse, then
a value close to the global maximum in $[0,T]$ will be found.
The following lemma explains why this is the case.

\begin{lemma}\label{lem:zetastaysbig}
Suppose
$|\zeta(\frac12+\i t_0)|=m_T:=\max_{t\in [0,T]} |\zeta(\frac12+\i t)|$.
There is an absolute constant $A>0$ such that if
$|t-t_0|<A/\log T$, then $|\zeta(\frac12+\i t)| > \frac12 m_T$.
\end{lemma}

\begin{proof}
We can
estimate the size of the derivative of the zeta-function near
$\frac12+\i t_0$ using Cauchy's theorem.
If we integrate around a circle of size $1/\log T$ and use the
functional equation, we find
that there exists an absolute constant $c_1$ such that if
$|s-(\frac12+\i t_0)|<c_1/\log T$, then $\zeta^{\prime}(s) \ll m_T
\log T$. This gives the lemma.
\end{proof}

As a random matrix model for  $|Z_X(\tfrac12+\i t)|$
when it is sampled at evenly spaced points, one might
consider the largest value of $K=K(N, X)$ such that
\begin{equation}
\Prob\left\{ \max_{1\leq j \leq N^c \exp(N\log X)}
|\Lambda_{U_j}(0)| \leq K \right\} \to 1.
\end{equation}
The following theorem determines $K$ explicitly as a function of $N$
and  $X$  and shows that
such sampling is sufficient to capture the large values.
\begin{theorem}\label{thm:max size Lambda}
Let $M=N^c e^{N\log X}$, where $2<X<N$ and $c>0$ is fixed.
If
\begin{equation}\label{eq:K2}
K_\varepsilon(N)=\exp\left(
(\frac{1}{\sqrt{\mathstrut 2}}+\varepsilon) \sqrt{\mathstrut N\log N\log X}
\right),
\end{equation}
and $U_1,\dots,U_M$ are chosen independently from $U(N)$, then as
$N\to\infty$,
\begin{equation}
\Prob\left\{ \max_{1\leq j \leq M} |\Lambda_{U_j}(0)| \leq
K_\varepsilon(N) \right\} \to 1
\end{equation}
for all $\varepsilon>0$ and for no $\varepsilon<0$.
\end{theorem}

\begin{proof}
As in the proof of Theorem~\ref{thm:max size max Lambda}, since the $U_j$ are independent we have
\begin{equation}\label{eq:large values Lambda(0)}
\log \Prob\left\{ \max_{1\leq j \leq M} |\Lambda_{U_j}(0)| \leq
K_\varepsilon(N) \right\} = M
\log\left(1-\Prob\left\{|\Lambda_{U}(0)| > K_\varepsilon(N)
\right\} \right) .
\end{equation}
Theorem 3.5 of \cite{HKO2} asserts that if $\delta>0$ is fixed  and
$\exp(N^\delta) \leq K \leq \exp(N^{1-\delta})$, then
\begin{equation}\label{eq:LDP Lambda}
\Prob\left\{ |\Lambda_{U}(0)| > K \right\} =
\exp\left(-\frac{\log^2 K}{ \log N - \log\log K}
(1+o_\delta(1))\right)
\end{equation}
as $N\to\infty$. Hence,  if $M=N^c \exp(N \log X)$
and $K_\varepsilon(N)$ is given by \eqref{eq:K2}, then the
left-hand side  of \eqref{eq:large
values Lambda(0)} tends to zero for all $\varepsilon>0$, but for
no $\varepsilon<0$.
\end{proof}

Note that the statements of Theorems~\ref{thm:max size max Lambda}
and~\ref{thm:max size Lambda} are almost identical and, in
particular, one can capture the largest values of $|\Lambda_U|$,
hence $|Z_X|$, just by sampling at individual points; it is  not
necessary to find the maxima of the individual polynomials.
This is significant for our modeling of the prime
contribution~$P_X$, for in that case we are only able to sample at
individual points, and there is nothing comparable to a sequence
of polynomials over which we can maximize individually.

%%%%%%%%%%%%%%%%%%%%%%

\subsection{Probabilistic model  for large values of $P_X$}\label{sect:prime model}

The material in this section was provided to us by Granville and Soundararajan.

First note that
\begin{eqnarray}
\log P_X(\tfrac12+\i t)
&=&\sum_{p\le X} \frac{1}{p^{\frac12+\i t}}
    + O\left(\sum_{p\le \sqrt{\mathstrut X}} \frac{1}{p}\right) \cr
&=&\sum_{p\le X} \frac{1}{p^{\frac12+\i t}} + O(\log\log X).
\end{eqnarray}
Hence
\begin{eqnarray}\label{eqn:Pstar}
P_X(\tfrac12+\i t) &=& \exp\left(\sum_{p\le X} \frac{1}{p^{\frac12+\i t}} \right)
    \times \exp(O(\log\log X)) \cr
&=& \exp\left(P^*_X(\tfrac12+\i t) \right)
        \times \exp(O(\log\log X)) ,
\end{eqnarray}
say.
We will see that this approximation is adequate as long as
$X=\exp(o(\sqrt{\mathstrut \log T\log\log T}))$.

The method is based  on treating the $p^{-\i t}$ as independent random
variables.  The large values of $P^*_X(\tfrac12+\i t)$ can then be
obtained from the following lemma, the proof of which involves
calculating the moments of the distribution.
\begin{lemma}  Let $\{z_j\}$ be a sequence of independent random variables distributed
uniformly on the unit circle and let $\{a_j\}$ be a sequence of bounded
real numbers such that for all $n \ge 3$,
\begin{equation}
\frac{1}{V_J^{\frac{n}{2}}} \sum_{1\le j\le J} a_j^n\to 0
\end{equation}
as $J\to\infty$, where
\begin{equation}
V_J:= \sum_{1\le j\le J} a_j^2 .
\end{equation}
Then, as $J\to\infty$, the distribution of
\begin{equation}
Y_J:= \Re \sum_{1\le j\le J} a_j z_j
\end{equation}
tends to a Gaussian with mean 0 and variance $\frac12 V_J$.
\end{lemma}

Applying the lemma to $P^*_X(\tfrac12+\i t)$ with $z_j=p_j^{-\i
t}$, where $p_j$ is the  $j$th prime, and $a_j=1/\sqrt{\mathstrut
\mathstrut p_j}$, we see that $V_J\sim \log\log X$. For
$X=\exp(\sqrt{\mathstrut \log T})$ we model the maximum of
$|P^*_X(\tfrac12+\i t)|$ by independently choosing $T \log^c T$
values of $t$. By \eqref{eqn:Pstar} this yields
\begin{equation}\label{eq: max P}
\max_{t\in [0,T]} |P_X(\tfrac12+\i t)| = \exp\left(
(1+o(1))\sqrt{\mathstrut \tfrac12 \log T\log \log T}\right) .
\end{equation}
If
$X=\exp(\sqrt{\mathstrut \log T})$, the method of
Section~\ref{sect:char poly} predicts that
\begin{equation}\label{eq: max Z}
\max_{t\in [0,T]} |Z_X(\tfrac12+\i t)| =
O\left(
    \exp\left(\sqrt{\mathstrut \log T}\right)
  \right) .
\end{equation}
These estimates together with \eqref{eq:zeta= PZ} give an
alternative justification of Conjecture~\ref{conj:max value zeta}.

\subsection{Combining $Z_X$ and $P_X$}\label{sect:primes and zeros}

If $X=\exp(\log^\alpha T)$ with $0<\alpha<\tfrac12$, then the
largest values of $|Z_X|$ and $|P_X|$ are  approximately the same
size and both will contribute to the largest values of
$|\zeta(\tfrac12+\i t)|$.  Specifically, applying
Theorem~\ref{thm:max size max Lambda} with $N=\log T / \log X$ (so that
$Z_X$ is modeled by $\Lambda_U$) and $M=T\log X$ (so that we sample
enough characteristic polynomials to cover the critical line
between $t=0$ and $t=T$), the previous analysis using characteristic
polynomials predicts that $|Z_X(\tfrac12+\i t)|$ gets as large as
\begin{equation}
\exp\left(\frac{1}{\sqrt{\mathstrut 2}} \sqrt{\mathstrut (1-2\alpha)\log T\log\log
T}\right),
\end{equation}
and $|P_X(\tfrac12+\i t)|$ gets as large as
\begin{equation}
\exp\left(\sqrt{\mathstrut \alpha\log T\log\log T}\right).
\end{equation}
The product of these two quantities is larger than our conjectured
maximum of $|\zeta(\tfrac12 + \i t)|$, as it should be, because we do
not expect $|Z_X|$ and $|P_X|$ to attain their maximum values
simultaneously. Instead of multiplying the maxima, we must find the distribution
of the large values of the product $|Z_XP_X|$ in order to check
that our method is consistent throughout the range
$0<\alpha<\tfrac12$.   This is a calculation involving
the  tails of the  distributions of $Z_X$ and~$P_X$.

In Section~\ref{sect:prime model} we saw that, conjecturally, if
$X=\exp(\log^\alpha T)$, then the tail of  $\log P_X$ has
the same  distribution as a Gaussian with mean 0 and
variance~$\sigma^2_P=\frac12 \alpha \log\log T$.  For $Z_X$,
Lemma~\ref{lemma:maxthetadist} states that if $\delta>0$ is fixed
and $\delta\le \lambda \le 1-\delta$, then
\begin{equation}
\Prob\left\{ \max_\theta |\Lambda_{U}(\theta)| \geq
\exp(N^\lambda) \right\} = \exp\left(-\frac{N^{2\lambda}}{
(1-\lambda) \log N} (1+o(1))\right) .
\end{equation}
Since $\Lambda_{U}$ models $Z_X$, the lemma provides the tail of
the distribution of $\log|Z_X(\tfrac12+\i t)|$. Assuming that
$ Z_X $ and $ P_X $ are essentially statistically independent,
we should expect the distribution
of  $\log|Z_X| + \log|P_X|$ to be the convolution of the two
distributions. Hence, for the tail we  convolve the
tails of the two distributions. Thus, for large $K$ we expect that
\begin{eqnarray}\label{eqn:convolution}
\lefteqn{\frac{1}{T} \meas\left\{ 0 < t < T \ : \
    \log |P_X(\tfrac{1}{2}+\i t)| + \log|Z_X(\tfrac{1}{2}+\i t)|\geq
\log K \right\} } \cr
 &=& \int_{-\infty}^\infty
\exp\left(-\(\frac{(\log K-x)^2}{\alpha\log\log T} +
    \frac{x^2}{(1-\alpha)\log\log T - \log x}\)
(1+o(1))\right) \d x \cr
 & =& \int_{-\infty}^\infty
\exp\left(- f_K(x) (1+o(1))\right) \d x ,
\end{eqnarray}
say.
By the saddle point method, this equals
\begin{equation}
\exp\left(- f_K(x_0)
(1+o(1))\right) ,
\end{equation}
where $x_0$ is such that $f_K'(x_0)=0$.
If we
write $K = \exp\left(d\sqrt{\mathstrut \log T \log\log
T}\right)$, then solving
\begin{multline}
0=f_K'(x_0) = \frac{2\left(x_0-d\sqrt{\mathstrut \log T \log\log
T}\right)}{\alpha \log\log T} +
\frac{2x_0}{(1-\alpha)\log\log T - \log x_0} \\
+ \frac{x_0}{((1-\alpha)\log\log T - \log x_0)^2}
\end{multline}
yields the solution $x_0 \sim d (1-2\alpha)\sqrt{\mathstrut \log T \log\log
T}$ (which is justified so long as $0<\alpha<1/2$).
Thus,
\begin{equation}
f_K(x_0) =  (2+o(1))\, d^2\log T  ,
\end{equation}
and this leads to the following conjecture:
\begin{conj}
For $d>0$ fixed and $T\to\infty$, we have
\begin{equation}
\frac{1}{T} \meas \left\{ 0<t<T  :  |\zeta(\tfrac12+\i t)| >
\exp(d \sqrt{\log T \log\log T}) \right\}   =    \exp\left( -2d^2
\log T (1 + o(1)) \right) .
\end{equation}
\end{conj}

By Lemma~\ref{lem:zetastaysbig}, $|\zeta(\frac12+\i t)|$ is close
to its maximum value over a window of size $C/\log T$, so we wish
to find the smallest $d$ such that
\begin{equation}
\meas\left\{ 0<t<T \ : \  |\zeta(\tfrac{1}{2}+\i t)| \geq \exp(d
\sqrt{\log T \log\log T}) \right\} \ll  \frac{1}{\log T},
\end{equation}
that is, such that  $T \log T \exp(-2d^2 (1+o(1))\log T) \ll  1$.
This happens if $d=\sqrt{\mathstrut \frac12}+\varepsilon$ for any
$\varepsilon>0$, but for no $\varepsilon<0$. Once more this gives
Conjecture~\ref{conj:max value zeta}, and this time in a way that
is independent of the choice of~$X=\exp(\log^{\alpha}T)$ for
$0<\alpha<1/2$.

%%%%%%%%%%%  !!!!!!!!!!!  %%%%%%%%%%%%%%%
%%%%%%%%%%%  New Section  %%%%%%%%%%%%%%%
%%%%%%%%%%%  !!!!!!!!!!!  %%%%%%%%%%%%%%%

\section{Bounds based on conjectures for
moments}\label{sect:moments}

In this section we obtain Conjecture~\ref{conj:max value zeta}
by using conjectures for moments of the
zeta-function.  Our method also leads to limits on the
possible uniformity of the conjectured moments.

Our approach here is based on the work of
Conrey and Gonek~\cite{ConGon}.
Let
\begin{equation}
m_T := \max_{t \in [0,T]} |\zeta(\tfrac12+\i t)|,
\end{equation}
and note that we have the trivial inequality
\begin{equation}\label{eq:trivial upper bound}
m_T^{2k} \geq \left( \frac{1}{T} \int_0^T |\zeta(\tfrac12 +\i
t)|^{2k} \;\d t \right).
\end{equation}
It follows that estimates for the right-hand side of the
inequality imply lower bounds for the maximum size of the
zeta-function.

Keating and Snaith~\cite{KeaSna}
used random matrix theory to conjecture that if
$k>-1/2$ is fixed, then as $T\to\infty$,
\begin{equation}
 \label{conj:KS}
\frac{1}{T} \int_0^T |\zeta(\tfrac12+\i t)|^{2k} \;\d t \sim
\frac{G^2(k+1)}{G(2k+1)} \, a(k) \log^{k^2} T,
\end{equation}
where $G$ is the Barnes $G$-function, and
\begin{equation}
a(k) = \prod_{\substack{p\\ \text{prime}}}
\left(1-\frac{1}{p}\right)^{k^2} \sum_{m=0}^\infty
\left(\frac{\Gamma(m+k)}{m! \ \Gamma(k)}\right)^{\!2} p^{-m} \; .
\end{equation}
This conjecture is for $k$ fixed, but we would like to let
$k\to\infty$, because the $2k$th root of the right-hand side of
\eqref{eq:trivial upper bound} then actually tends to $m_T$.
Thus, we would like to know how large $k$ can be
as a function of $T$.

Conrey and Gonek showed that if formula~(\ref{conj:KS}) holds
for $k=\sqrt{\mathstrut \log T/ \log\log T}$ then
\begin{equation}\label{eq:ConGon max1}
m_T \geq \exp\left(C_1 \sqrt{\mathstrut \frac{\log
T}{\log\log T}}\right),
\end{equation}
and if it holds for $k$ as large as $\log T/ \log\log T$
then
\begin{equation}\label{eq:ConGon max}
m_T \geq \exp\left(C_2 \frac{\log
T}{\log\log T}\right) ,
\end{equation}
where $C_1$ and $C_2$ are given explicitly.
Hughes \cite{Hug} gave a convexity argument to show that
formula~(\ref{conj:KS}) must fail before $k=\log T/ \log\log T$.
However, using the last point at which convexity holds
for~(\ref{conj:KS}), one still obtains \eqref{eq:ConGon max},
but with a smaller constant~$C_2$.

In all of these cases one only requires  a lower bound for
the right-hand side of~(\ref{conj:KS}).
Our approach here is to use the mean value formula \eqref{conj:KS} to
obtain upper bounds instead of lower bounds for $m_T$.
As a consequence, we  also obtain restrictions on the
possible range of validity of~\eqref{conj:KS} for $k$ growing with $T$.
Specifically we prove the following:
\begin{theorem} \label{thm:KS fails}
Formula \eqref{conj:KS} does not hold for
\begin{equation}
k\geq (2\sqrt{\mathstrut 2}+\varepsilon) \sqrt{\mathstrut
\frac{\log T}{\log\log T}}
\end{equation}
for any fixed $\varepsilon >0$.
\end{theorem}

Our method allows us to get upper bounds for $m_T$ and,
in particular,  we obtain
\begin{theorem} \label{thm:max size zeta}
If formula \eqref{conj:KS} holds for
$k=\log^\delta T$ for some $\delta<\frac12$, then
\begin{equation}
m_T \ll \exp\left( \log^{1-\delta} T \right) .
\end{equation}
Moreover, if formula \eqref{conj:KS} holds for
$k=\sqrt{\mathstrut 2\log T/\log\log T}$, then
\begin{equation}
m_T \ll \exp\left(  \sqrt{\mathstrut \tfrac12 \log T \log\log T} +
O\left(\frac{\sqrt{\mathstrut \log T} \log\log\log T}{\sqrt{\mathstrut \log\log
T}}\right)\right) ,
\end{equation}
and if formula \eqref{conj:KS} holds for $k=\sqrt{\mathstrut 8\log
T/\log\log T}$, then
\begin{equation}
m_T \gg \exp\left( \sqrt{\mathstrut \tfrac12 \log T \log\log T} +
O\left(\frac{\sqrt{\mathstrut \log T} \log\log\log T}{\sqrt{\mathstrut \log\log
T}}\right)\right) .
\end{equation}
\end{theorem}

Theorem~\ref{thm:max size zeta} says that if formula
\eqref{conj:KS} holds true until $k = \sqrt{\mathstrut 8\log T / \log\log T}$
(after which we know it must fail), then we have
\begin{equation}
\label{eqn:trueorder}
m_T = \exp\left( \sqrt{\mathstrut \tfrac12 \log T \log\log T}
\left(1+O\left(\frac{\log\log\log T}{\log\log T}\right) \right) \right)
\end{equation}
Note that this implies the conjecture found in the previous section.
This is not surprising because, as we will show, the arithmetic factor
in the conjectured moments is smaller than the other factors.  Thus,
this bound is coming just from the random matrix model.

If the true order of the zeta-function is larger than the
bound in \eqref{eqn:trueorder}, then one would like to know
where our calculation fails.  Since
formula \eqref{conj:KS} is only the leading order term in the
asymptotic expansion for the $2k$th moment of the zeta-function,
it is possible that the lower order terms dominate when $k=\log^\delta T$.
However, this is unlikely.  In the random matrix case, the $2k$th
moment is given by
\begin{align}
\E\left\{\left|\Lambda_U(\theta\right|^{2k} \right\} &=
\frac{G^2(k+1)}{G(2k+1)} \frac{G(N+1) G(N+2k+1)}{G^2(N+k+1)}\\
&=\frac{G^2(k+1)}{G(2k+1)} \exp\left(k^2\log N + \frac{k^3}{N} -
\frac{7k^4}{12 N^2} + \frac{k^5}{2N^3} + \dots \right),
\end{align}
and one sees that the first term dominates even
for~$k$ as large as $N^\delta T$ provided that $\delta<1$.

In the zeta-function case, the complete main term of the $2k$th
moment has been conjectured (see~\cite{CFKRS}, Conjecture~1.5.1). One
can check that the contribution from the primes is bounded by
$\exp(c k^2)$, which is insufficient to affect the estimate for
$m_T$.  Thus, if our conjecture for the growth of $|\zeta(\frac12 +
\i t)|$ is incorrect, then the main term in the mean value must
take a new form for $k=\log^\delta T$ for some $\delta<\frac12$.
If~\eqref{eq:zeta_bigO} is the true maximal size, then by equation
\eqref{eq:320} the conjectured mean value can only hold for $k \ll
\log\log T$.

Our main tool for finding upper bounds is the following lemma.

\begin{lemma}\label{lem:lower bound}
For all positive real $k$ we have
\begin{equation}
m_T^{2k} \ll 2^{2k} \log T \int_0^T |\zeta(\tfrac12+\i t)|^{2k} \d t ,
\end{equation}
where the implied constant is absolute.
\end{lemma}

\begin{proof}

Suppose that  $|\zeta(\frac12+\i t_0)|=m_T$, where $0<t_0<T$.
By Lemma~\ref{lem:zetastaysbig}
there is an absolute constant $A>0$ such that if $|t - t_0|\le
A/\log T$ then $|\zeta(\frac12+\i t)| \ge \frac12 m_T$. This gives
\begin{eqnarray}
\int_0^T |\zeta(\tfrac12+\i t)|^{2k} \d t
&\ge &
\int_{t_0-A/\log T}^{t_0+A/\log T} |\zeta(\tfrac12+\i t)|^{2k} \d t \cr
&\ge &
\frac{A}{\log T} \left(\frac{m_T }{2}\right)^{2k} ,
\end{eqnarray}
as claimed.
\end{proof}

\begin{proof}[Proof of Theorems~\ref{thm:KS fails} and~\ref{thm:max size zeta}]
By Lemma~\ref{lem:lower bound} and \eqref{eq:trivial upper bound}
there exists an absolute constant $C$ such that
\begin{equation}
\label{eqn:doubleinequality}
\left(\frac{1}{T} \int_0^T |\zeta(\tfrac12+\i t)|^{2k} \;\d
t\right)^{1/2k} \leq m_T \leq
2(CT \log T)^{1/2\ell}  \left(\frac{1}{T} \int_0^T
|\zeta(\tfrac12+\i t)|^{2\ell} \;\d t \right)^{1/2\ell} .
\end{equation}
We use these inequalities to prove Theorem~\ref{thm:max size zeta} first.

The asymptotic expansion for the Barnes $G$-function (see Barnes
\cite{Bar})  gives
\begin{equation}\label{eq:G asy}
\frac{G(1+k)^2}{G(1+2k)} = \exp\left(k^2\left(-\log k +
\frac{3}{2} - 2\log 2\right) - \frac{1}{12} \log k  + \frac{1}{12}
\log 2 + \zeta'(-1) + O\left(\frac{1}{k}\right) \right)
\end{equation}
for $k\geq 1$.
Furthermore, Conrey and Gonek \cite{ConGon} have shown that
\begin{equation}\label{eq:a asy}
\log a(k) \sim -k^2\log(2e^\gamma\log k) + o(k^2) \textrm{ for } k
\rightarrow \infty .
\end{equation}
Thus, if  \eqref{conj:KS} holds, then
\begin{equation}
\label{eqn:logKs}
\log  \left(\frac{1}{T} \int_0^T
|\zeta(\tfrac12+\i t)|^{2\ell} \;\d t \right)^{1/2\ell}
=
\frac12 \ell \log\log T
-\frac{1}{2} \ell \log \ell +O( \ell \log\log \ell) .
\end{equation}
It therefore follows from~\eqref{eqn:doubleinequality} that
\begin{equation}\label{eq:320}
\log m_T \ll
\frac{\log T+\log\log T}{2\ell} +\frac12 \ell \log\log T
-\frac{1}{2} \ell \log \ell +O( \ell \log\log \ell) .
\end{equation}
Setting $\ell=\log^\delta T$, we find that
\begin{equation}
m_T \ll \exp\left(\frac{1}{2} \log^{1-\delta} T +O(\log^\delta T\log\log T)
\right),
\end{equation}
which gives the first inequality in Theorem~\ref{thm:max size zeta}.
Setting $\ell=c \sqrt{\frac{\mathstrut \log T}
{ \mathstrut \log\log T}}$, we find that
\begin{equation}\label{eqn:secondestimate}
m_T \ll \exp\left(\left(\frac{1}{2c} +
\frac{c}{4}\right)\sqrt{\mathstrut \log T \log\log T} +
O\left(\frac{\sqrt{\mathstrut \log T} \log\log\log T}{\sqrt{\mathstrut \log\log
T}}\right)\right),
\end{equation}
which is minimized by taking $c=\sqrt{\mathstrut 2}$.  This gives the second estimate
in Theorem~\ref{thm:max size zeta}.

By~\eqref{eqn:doubleinequality} and
\eqref{eqn:logKs}
we also have
\begin{equation}
m_T\gg \exp\left(\frac12 k \log\log T - \frac12 k \log k + O(k\log\log k) \right).
\end{equation}
If we set $k=c \sqrt{\frac{\mathstrut \log T}
{ \mathstrut \log\log T}}$,
we find that
\begin{equation}
m_T\gg \exp\left( \frac{c}{4}\sqrt{\mathstrut \log T \log\log T}
+ O\left(\frac{\sqrt{\mathstrut \log T} \log\log\log T}{\sqrt{\mathstrut
\log\log T}}\right)\right).
\end{equation}
Choosing $c=2\sqrt{\mathstrut 2}$, we obtain the third inequality
 in Theorem~\ref{thm:max size zeta}.
If $c>2\sqrt{\mathstrut 2}$,  this contradicts~\eqref{eqn:secondestimate}
and  thereby establishes
Theorem~\ref{thm:KS fails}.
\end{proof}

\section{Bounds for $S(t)$}\label{sect:S}

Recall that $S(t)$ is the error term in the counting function for
the number of non-trivial zeros of the zeta-function with imaginary part less
than $t$. It may also be expressed as
$S(t) = \frac{1}{\pi} \Im\log \zeta(\tfrac12+\i t)$ (see Titchmarsh \cite{Ti}).
Since $ \zeta(\frac12+\i t)$ is essentially $P_X(\frac12+\i t) Z_X(\frac12+\i t)$
and $Z_X$ is modeled by $ \Lambda_U(\theta)$,
one would expect that if $X$  is  sufficiently small, so that the contribution
from $P_X$ is negligible, then
 $\frac{1}{\pi} \Im\log \zeta(\tfrac12+\i t)$ too can
be modeled by random matrix theory, in particular, by
\begin{equation}
\frac{1}{\pi} \Im \log \Lambda_U(0).
\end{equation}
Evidence for this is presented in \cite{KeaSna}.  This is the basis for our
Conjecture~\ref{conj:max value S}.

\begin{theorem}\label{thm:im log Lambda}
Set
\begin{equation}\label{eq:K for im log}
K_\varepsilon(N)=\left(\frac{1}{\sqrt{\mathstrut 2}}+\varepsilon\right)
\sqrt{\mathstrut N\log N}
\end{equation}
and $M=N^c e^N$, where $c>0$ is fixed. If $U_1,\dots,U_M$ are
chosen independently from $U(N)$, then as $N\to\infty$,
\begin{equation}
\Prob\left\{ \max_{1\leq j \leq M} \Im \log \Lambda_{U_j}(0) \leq
K_\varepsilon(N) \right\} \to 1
\end{equation}
for all $\varepsilon>0$ and for no $\varepsilon<0$.
\end{theorem}

\begin{proof}
This follows along the same lines as the proof of
Theorem~\ref{thm:max size Lambda}. Since we are making independent
choices,
\begin{equation}
\Prob\left\{ \max_{1\leq j \leq M} \Im\log\Lambda_{U_j}(0) \leq K
\right\} = \Prob\left\{ \Im\log\Lambda_{U}(0) \leq K \right\}^M .
\end{equation}
For this to tend to $1$, we need
\begin{equation}\label{eq:Im log condition on max}
M \log \Prob\left\{ \Im\log\Lambda_{U}(0) \leq K \right\} \to 0 .
\end{equation}
By Theorem 3.6 of Hughes, Keating and O'Connell \cite{HKO2},
if $K=N^\lambda$, where $\delta<\lambda<1-\delta$ and
$\delta>0$ is fixed, then
\begin{equation}
\Prob\left\{ \Im\log\Lambda_{U}(0) \leq K) \right\}= 1 -
\exp\left(-\frac{K^2}{ \log N - \log K} \left(1+o_\delta(1)\right)
\right) .
\end{equation}
One can easily check that if $K_\varepsilon(N)$  is given by
\eqref{eq:K for im log}, then \eqref{eq:Im log condition on max}
holds for all $\varepsilon>0$, but for no $\varepsilon<0$.
\end{proof}

Conjecture~\ref{conj:max value S} now follows in the same manner as
the justification of Conjecture~\ref{conj:max value zeta}; that is, by controlling
the prime contribution from $\Im\log  P_X(\tfrac12+\i t)$.

\section{Other families and other arguments}\label{sect:other}

\subsection{Other families: symplectic and orthogonal}

The analogue of Gonek,  Hughes,  and Keating's approximation to
the zeta-function has not yet been extended to the case of other
$L$-functions near the critical point.  However, it is believed
that the characteristic polynomials of symplectic (or orthogonal)
matrices model the central value of $L$-functions taken from a
symplectic (or orthogonal) family of $L$-functions \cite{KS2}.
Thus, the methods developed in section~\ref{sect:char poly} can be
applied. Moreover, we can still estimate the maximal size of
critical values by using a partial Euler product and modifying the
method of Granville and Soundararajan. Finally, we can also apply the
method involving mean values.

We give as an example finding the large values of the characteristic
polynomials of the symplectic group at the critical point. The orthogonal
family is treated in an almost identical way. The characteristic
polynomial of an $N\times N$ symplectic matrix ($N$ must be even)
with eigenvalues $e^{\pm\i\theta_n}$ is
\begin{equation}
Z(U,0) = \prod_{j=1}^{N/2} (1-e^{\i\theta_n})(1-e^{-\i\theta_n}).
\end{equation}
Keating and Snaith \cite{KS2} calculated the moment generating
function and found that
\begin{equation}
\E_{Sp(N)}\left\{Z(U,0)^s \right\} = 2^{Ns} \prod_{j=1}^{N/2}
\frac{\Gamma(N/2+j+1) \Gamma(s+j+1/2)}{\Gamma(j+1/2)
\Gamma(s+N/2+j+1)}.
\end{equation}
A long but straightforward calculation using Stirling's asymptotic
series for the gamma function shows that if $\delta>0$ is fixed
and $\delta<\lambda<1-\delta$, then for $A(N) = N^{\lambda}$,
$B(N) = \frac{N^{2\lambda}}{(1-\lambda)\log N}$, and fixed $s \geq
0$, we have
\begin{equation}\label{eq:symplectic log gen func}
\lim_{N\to\infty} \frac{1}{B} \log \E_{Sp(N)}\left\{Z(U,0)^{s B /
A} \right\} = \frac{1}{2} s^2.
\end{equation}
{}From this, large deviation theory (for example, see \cite{DZ})
allows us to deduce that if
$\exp(N^\delta) \leq K \leq \exp(N^{1-\delta})$, then
\begin{equation}\label{eq:symplectic rate func}
\Prob_{Sp(N)}\left\{ Z(U,0) > K \right\} = \exp\left(-\frac{\log^2
K}{ 2\log N -2\log\log K} (1+o_\delta(1))\right)
\end{equation}
as $N\to\infty$. Comparing this with \eqref{eq:LDP Lambda}, the
analogous statement for the  unitary group, we note the extra
factor of $2$ in the denominator. This difference explains why the
constant $B$ in Conjecture~\ref{conj:max value L} equals $1$
rather than $1/2$.

We now see,  by  methods identical to those of the previous
section, that if $M=N^c e^{N}$ for any fixed $c\geq 0$, and if
$K_\varepsilon=\exp\left((1+\varepsilon)\sqrt{\mathstrut N \log
N}\right)$, then if $U_1,\dots, U_M$ are chosen independently from
$Sp(N)$,
\begin{equation}
\Prob_{Sp(N)}\left\{ \max_{1\leq j \leq M} Z(U_j,0) \leq
K_\varepsilon(N) \right\} \to 1
\end{equation}
as $N\to\infty$ for all $\varepsilon>0$ and for no
$\varepsilon<0$.

Consider for instance the family of all quadratic Dirichlet
$L$-functions $L(s,\chi_d)$. For characters with modulus around
$D$, random matrix theory suggests (see Keating and Snaith \cite{KS2})
that $N = \log D$ is the correct identification between the size of the
matrix and the conductor (though note that in \cite{KS2} $N$ is \emph{half}
the size of the symplectic matrix). Furthermore, it is well known that there
are asymptotically $6 D / \pi^2$ primitive discriminants less than $D$.
Thus, we conjecture that
\begin{equation}
\max_{\substack{|d| \leq D \\ \chi_d \text{ real}}}
|L(\tfrac12,\chi_d)| = \exp\left(  (1+o(1)) \sqrt{\mathstrut \log
D \log\log D}\right) .
\end{equation}

Similarly,  the moment generating function
has been calculated for the orthogonal case (see \cite{KS2})
and, if $N$ is even, we have
\begin{equation}
\E_{SO(N)}\left\{Z(U,0)^s \right\} = 2^{Ns} \prod_{j=1}^{N/2}
\frac{\Gamma(N/2+j-1) \Gamma(s+j-1/2)}{\Gamma(j-1/2)
\Gamma(s+N/2+j-1)}.
\end{equation}
Equations \eqref{eq:symplectic log gen func} and
\eqref{eq:symplectic rate func} apply to the orthogonal case
without change, so by the same reasoning as previously, if $M=N^c
e^{N}$ for any fixed $c\geq 0$, and
$K_\varepsilon=\exp\left((1+\varepsilon)\sqrt{\mathstrut N \log N}\right)$, then
\begin{equation}
\Prob_{SO(N)}\left\{ \max_{1\leq j \leq M} Z(U_j,0) \leq
K_\varepsilon(N) \right\} \to 1
\end{equation}
as $N\to\infty$ for all $\varepsilon>0$, but for no
$\varepsilon<0$.

Next we consider how to adapt the Granville-Soundararajan argument
involving the product over primes
to the symplectic case. We require the following lemma.
\begin{lemma}
Let $\{x_j\}$ be a sequence of independent \emph{real} random variables
with mean~0 and variance~1, and let $\{a_j\}$ be a bounded sequence
of real numbers such that for all $n \ge 3$,
\begin{equation}
\frac{1}{V_J^{\frac{n}{2}}} \sum_{1\le j\le J} a_j^n\to 0
\end{equation}
as $J\to\infty$, where
\begin{equation}
V_J:= \sum_{1\le j\le J} a_j^2 .
\end{equation}
Then as $J\to\infty$, the distribution of
\begin{equation}
Y_J:= \sum_{1\le j\le J} a_j x_j
\end{equation}
tends to a Gaussian with mean 0 and variance $V_J$.
\end{lemma}

Just as in the treatment involving characteristic polynomials, the
fact that the variance is $V_j$ for these families instead of
$V_j/2$ leads to the constant $B=1$ in Conjecture~\ref{conj:max
value L}, instead of $B=1/2$ for the unitary family dealt with
previously.

\subsection{Other arguments, for}

Although the conjectures in this paper are based on very recent
work,
Hugh Montgomery has pointed out to us that a similar conjecture can be obtained
by viewing $\log|\zeta(\frac12+\i t)|$ as a Gaussian distributed random
variable with variance $C \log\log T$, where one estimates $m_T$
by sampling $T^A$ times.

Soundararajan suggests a different way to use the moments of the
zeta-function to conjecture an upper bound.
The proof of Lemma~\ref{lem:lower bound} showed that if
\begin{equation}
\frac{1}{T}\meas\left\{ t \in [0,T] \ : \ |\zeta(\tfrac12+\i t)|
\geq \tau \right\} \leq \frac{1}{T},
\end{equation}
then $m_T \leq 2\tau$.  For when $|\zeta(\tfrac12+\i t)|$ is very
large, it must remain large over an interval of size $c/\log T$.
Now
\begin{equation}
\tau^{2k} \meas\left\{ t \in [0,T] \ : \ |\zeta(\tfrac12+\i t)|
\geq \tau \right\} \leq \int_0^T |\zeta(\tfrac12 +\i t)|^{2k} \;\d t \,,
\end{equation}
so if equation \eqref{conj:KS} holds, then we have
\begin{equation}
\frac{1}{T}\meas\left\{ t \in [0,T] \ : \ |\zeta(\tfrac12+\i t)|
\geq \tau \right\} \leq \tau^{-2k} \frac{G^2(k+1)}{G(2k+1)} (\log
T)^{k^2}.
\end{equation}
The right-hand side is less than $c/(T\log T)$  (which means  there
is only one place where the maximum occurs) when
\begin{equation}
\tau \geq \exp\left(\frac{\log T}{2k} + \frac{k}{2} \log\log T -
\frac{k}{2} \log k\right) .
\end{equation}
The minimum of this is $\exp\left(\sqrt{\mathstrut \frac12\log T \log\log
T}\right)$  and it occurs when $k=\sqrt{\mathstrut 2 \log T / \log\log T}$.

\subsection{Other arguments, against}

We now discuss  potential arguments against the conjectures
in this paper. One possibility,  so fundamental that
it cannot be addressed, is that the
large values of an $L$-function may be so rare that
these statistical models cannot detect them. Indeed, since these are
problems in number theory, there may be number-theoretic
constructions of large values which contradict our conjectures.
The two examples below,  due to Brian Conrey, suggest the kinds
of things we have in mind.

The first argument invokes an
analogy with the divisor function $d(n)=\sum_{d|n} 1$ and the
related function $\omega(n)=\sum_{p|n} 1$, where here the sum is
over  prime divisors of $n$. Since $\omega(n)$ is $\log\log n$ on average and
has a Gaussian distribution,  the relation
$d(n) = 2^{\omega(n)}$ for square-free $n$ might lead one to conjecture that
for such $n$, $d(n)$ is bounded by
\begin{equation}
\exp(c \sqrt{\mathstrut  \log n/\log\log n  } ).
\end{equation}
However,  we know how to
construct large values of $d(n)$, and it is easy to see that
for $n$ squarefree, $d(n)$ can get as large as
\begin{equation}
\exp(c \log n/\log\log n   ).
\end{equation}

The second argument  concerns the
Fourier coefficients $a_n$ of cusp forms.  For integer weight cusp
forms, rescaled so that for $p$ prime we have $|a_p|\le 2$, the
coefficients $a_n$ can get as large as
\begin{equation}
\exp(c \log n/\log\log n   ).
\end{equation}
In other words, they can get about as large as~$d(n)$. The
question is: can the coefficients of half-integral weight forms
also get this large?  If they can, then our conjecture on the
maximal size of the critical values of a symplectic family of
$L$-functions is incorrect.  For if $f\in
S_k(\Gamma_0(N))$, then $L_f(\frac12,\chi_d) = c_d^2/\sqrt{\mathstrut d}$,
where  $c_d$ is a Fourier coefficient of the  half-integral weight
form associated with $f$ by the Shimura correspondence.

Our methods cannot be directly applied to produce a version of
Conjecture~\ref{conj:max value S} for symplectic and orthogonal
families. But if one were to believe that for a family $\mathcal
F$ of $L$-functions with $c(F)$ denoting the conductor of~$F\in
\mathcal F$,
\begin{equation}\label{eq:gen_conj_B}
\limsup_{c(F)\to\infty} \frac{\Im \log
F(\tfrac12)}{\sqrt{\mathstrut \log c(F) \log\log c(F) }} =
\sqrt{\mathstrut B }
\end{equation}
where $B=1/2$ for unitary families and $B=1$ for symplectic and
orthogonal families, then since the rank of an elliptic curve is
related to the order of vanishing of its associated $L$-function,
this could lead to new information about large ranks. That is, if
\eqref{eq:gen_conj_B} is true, then it suggests that for rational
elliptic curves we have
\begin{equation}
\limsup_{c_E \to \infty} \frac{\rank(E)}{\sqrt{\mathstrut \log c_E
\log\log c_E}} = 1 ,
\end{equation}
where $c_E$ is the conductor of~$E$.
Note that this is smaller than the ranks of elliptic curves found
by Ulmer~\cite{Ul} in the function field case.
%This implies at
%least one of two things: Either $L$-functions over function fields
%are atypical of $L$-functions over number fields, or the
%generalization of Conjecture~\ref{conj:max value S} is not given
%by \eqref{eq:gen_conj_B} for non-unitary families.

%%%%%%%%%%%%%%%

\appendix
\section{The tail of the distribution of $\max_\theta |\Lambda_U(\theta)|$}
\label{sect:proofs}

Here we prove  the  random matrix polynomial result used in
Section~\ref{sect:new model approach}.

\begin{lemma}\label{lemma:maxthetadist}
If $\delta>0$ is fixed and $\delta\le \lambda \le 1-\delta $,
then
\begin{equation}\label{eqn:lambdatail}
\Prob\left\{ \max_\theta |\Lambda_{U}(\theta)| \geq
\exp(N^\lambda) \right\} = \exp\left(-\frac{N^{2\lambda}}{
(1-\lambda) \log N} (1+o(1))\right) .
\end{equation}
\end{lemma}

\begin{proof}[Proof of Lemma~\ref{lemma:maxthetadist}]

Bernstein's inequality for polynomials implies that for any matrix $U$,
\begin{equation}
\max_\theta | \Lambda'_U(\theta) | \leq N \max_\theta |
\Lambda_U(\theta) | .
\end{equation}
Thus, if $\phi$ is a point at which the maximum of
$|\Lambda_U(\theta)|$ occurs, and if we indicate the maximum
by $m_U$, then
for $|\theta-\phi| \leq 1/N$,
\begin{equation}
|\Lambda_U(\theta)| \geq m_U -|\theta-\phi| N m_U.
\end{equation}
It follows that
\begin{align}
\int_0^{2\pi} |\Lambda_U(\theta)|^{2k} \;\d \theta &\geq m_U^{2k}
\int_{-1/N}^{1/N} \left(1-|x| N\right)^{2k}\;\d x\\
&=m_U^{2k} \frac{2}{N} \frac{1}{2k+1}\, .
\end{align}
Combining this with the trivial lower bound for $m_U$, we find that
\begin{equation}\label{inequality}
\frac{1}{2\pi} \int_{0}^{2\pi} |\Lambda_U(\theta)|^{2k} \;\d
\theta \leq m_U^{2k} \leq \frac{2k+1}{2} N \int_{0}^{2\pi}
|\Lambda_U(\theta)|^{2k} \;\d \theta .
\end{equation}
This  bound holds for any matrix. We now average over all $N
\times N$ unitary matrices  with respect to Haar measure. That is,
we calculate the expectation $\E_N $ of
$|\Lambda_U(\theta)|^{2k}$. Set
\begin{equation}
\E_N \left\{ |\Lambda_U(\theta)|^{2k} \right\} = M_N(2k).
\end{equation}
Keating and Snaith~\cite{KeaSna} have shown that
\begin{equation}
M_N(2k) = \frac{G^2(k+1)}{G(2k+1)} \frac{G(1+N)
G(1+N+2k)}{G^2(1+N+k)},
\end{equation}
where $G$ is the Barnes $G$-function. Note that this is
independent of $\theta$. Therefore, by \eqref{inequality}
\begin{equation}
M_N(2k) \leq \E\left\{ m_U^{2k} \right\} \leq \pi(2k+1) N M_N(2k) .
\end{equation}
Hughes, Keating and O'Connell~\cite{HKO2} have shown that if $A(N) =
N^{\lambda}$ with $\delta<\lambda<1-\delta$ and $\delta>0$
fixed, and if
\begin{equation}
\label{eqn:BN} B(N) = \frac{N^{2\lambda}}{(1-\lambda)\log N} ,
\end{equation}
then for $s\geq 0$,
\begin{equation}
\lim_{N\to\infty} \frac{1}{B(N)} \log M_N\left(\frac{s B(N)
}{A(N)}\right) = \frac14 s^2 .
\end{equation}
Since
\begin{align}
\frac{1}{B(N)} \log M_N(s B(N) / A(N)) \leq  \mathstrut&\frac{1}{B(N)} \log
\E\left\{ m_U^{s B(N) / A(N)} \right\} \cr
\leq  \mathstrut&\frac{1}{B(N)} \log M_N(s B(N) / A(N)) + O\left(\frac{(\log
N)^2}{N^{2\lambda}}\right) ,
\end{align}
we   conclude that for $s\geq 0$,
\begin{equation}
\lim_{N\to\infty} \frac{1}{B(N)} \log \E\left\{ \exp\left(\frac{s
B(N) \log (\max_\theta |\Lambda_U(\theta)|)}{ A(N)}\right)
\right\} = \frac14s^2 .
\end{equation}
{}From this, large deviation theory (see, for example, \cite{DZ})
allows us to deduce that
\begin{equation}\label{eq:LDP for max}
\lim_{N\to\infty} \frac{1}{B(N)} \log \Prob\left\{ \log
\max_\theta |\Lambda_U(\theta)| \geq A(N) \right\} = -1 .
\end{equation}
Inserting  $A(N)=N^\lambda$ and $B(N)$ from (\ref{eqn:BN}), we
obtain the statement in the lemma.
\end{proof}


\begin{thebibliography}{99}
\bibitem{Bal} R. Balasubramanian, ``On the frequency of Titchmarsh's phenomenon for $\zeta(s)$ IV'', {\em Hardy-Ramanujan J.} \textbf{9} (1986) 1--10

\bibitem{BalRam} R. Balasubramanian and K. Ramachandra, ``On the frequency of Titchmarsh's phenomenon for $\zeta(s)$ III'', {\em Proc. Indian Acad. Sci. A} \textbf{86} (1977) 341--351

\bibitem{Bar} E.W. Barnes, ``The theory of the $G$--function'', \textit{Quart. J. Pure Appl. Math.} \textbf{31} (1899) 264--314

\bibitem{CF} J.B. Conrey, D.W. Farmer, ``Mean values of L-functions and symmetry'',
\textit{Internat. Math. Res. Not.} \textbf{17} (2000) 883--908

\bibitem{CFKRS} J.B. Conrey, D.W. Farmer, J. Keating, M. Rubinstein, and N.C. Snaith,
``Integral moments of L-functions'', \textit{Proc. London Math.
Soc. (3)} \textbf{91} (2005) 33--104

\bibitem{ConGon} J.B. Conrey and S.M. Gonek, ``High moments of the Riemann zeta function'', \textit{Duke Math. J.} \textbf{107} (2001) 577--604

\bibitem{DZ} A. Dembo and O. Zeitouni, \textit{Large Deviations Techniques and Applications, 2nd Ed.} (Springer-Verlag, 1998)

\bibitem{GHK} S.M. Gonek, C.P. Hughes and J.P. Keating, ``A Hybrid Euler-Hadamard product for the Riemann zeta function'', math.NT/0511182, to appear in Duke Math. J.

\bibitem{GS1}  A. Granville and K. Soundararajan, ``The distribution of
values of $L(1,\chi\sb d)$'',  \textit{Geom. Funct. Anal.}  13  (2003),
no. 5, 992--1028

\bibitem{GS2}  A. Granville and K. Soundararajan, ``Extreme
values of $|\zeta(1+\i t)|$'', math.NT/0501232.

\bibitem{Hug} C.P. Hughes, \textit{On the
Characteristic Polynomial of a Random Unitary Matrix and the
Riemann Zeta Function}, PhD Thesis (University of Bristol, 2001)

\bibitem{HKO2} C.P. Hughes, J.P. Keating and N. O'Connell, ``On the characteristic polynomial of a random unitary matrix'', \textit{Comm.. Math. Phys.} \textbf{220} (2001) 429--451

\bibitem{KeaSna} J.P. Keating and N.C. Snaith, ``Random matrix theory and $\zeta(1/2+\i t)$'', \textit{Comm. Math. Phys.} \textbf{214} (2000) 57--89

\bibitem{KS2} J.P. Keating and N.C. Snaith, ``Random matrix theory and $L$--functions at $s=1/2$'', \textit{Comm. Math. Phys.} \textbf{214} (2000) 91--110

\bibitem{K} T. Kotnik, ``Computational estimation of the order of $\zeta(\frac12+i t)$'', \textit{Math. Comp.} \textbf{73} no 246 (2004) 949-956

\bibitem{Mon} H.L. Montgomery, ``Extreme values of the Riemann zeta function'',
\textit{Comment. Math. Helvetici} \textbf{52} (1977) 511--518

\bibitem{MV} H.L. Montgomery and R.C. Vaughan, ``Extreme values of
Dirichlet $L$-functions at $1$'',  \textit{Number theory in progress,
Vol. 2 (Zakopane-Ko\'scielisko, 1997)},  1039--1052, de Gruyter, Berlin, 1999

\bibitem{S} Soundararajan, personal communication.

\bibitem{Ti} E.C. Titchmarsh, \textit{The Theory of the Riemann Zeta-Function} 2nd ed., revised by  D.R. Heath-Brown (Oxford Science Publications, 1986)

\bibitem{Ul} D. Ulmer, ``Elliptic curves with large rank over function fields'',
\textit{Ann. of Math.} \textbf{155} (2002), 295-315.


\end{thebibliography}
\end{document}